\input amstex
 \documentstyle{amsppt}
 \magnification=\magstep1
 \vsize=24.2true cm
 \hsize=15.3true cm
 \nopagenumbers\topskip=1truecm
 \headline={\tenrm\hfil\folio\hfil}

 \TagsOnRight

\hyphenation{auto-mor-phism auto-mor-phisms co-homo-log-i-cal co-homo-logy
co-homo-logous dual-izing pre-dual-izing geo-metric geo-metries geo-metry
half-space homeo-mor-phic homeo-mor-phism homeo-mor-phisms homo-log-i-cal
homo-logy homo-logous homo-mor-phism homo-mor-phisms hyper-plane hyper-planes
hyper-sur-face hyper-sur-faces idem-potent iso-mor-phism iso-mor-phisms
multi-plic-a-tion nil-potent poly-nomial priori rami-fication sin-gu-lar-ities
sub-vari-eties sub-vari-ety trans-form-a-tion trans-form-a-tions Castel-nuovo
Enri-ques Lo-ba-chev-sky Theo-rem Za-ni-chelli in-vo-lu-tion Na-ra-sim-han Bohr-Som-mer-feld}

\define\rest#1{_{\textstyle{|}#1}} 

\define\Span#1{\left<#1\right>} 
\define\na{\nabla}

\define\half{{\textstyle{1\over2}}}


\define\C{\Bbb C} 
\define\R{\Bbb R} 
\define\Z{\Bbb Z} 

\define\proj{\Bbb P} 

\define\sA{{\Cal A}} 
\define\sG{{\Cal G}} 
\define\sM{{\Cal M}} 

\define\al{\alpha}
\define\be{\beta}
\define\de{\delta}

\define\ga{\gamma}

\define\om{\omega}
\define\si{\sigma}
\define\De{\Delta}
\define\Ga{\Gamma}
\define\La{\Lambda}
\define\Om{\Omega}

\define\Si{\Sigma}





\define\W{\mathop{\bigwedge^2}} 



 \document

  \topmatter
  \title  Two lectures on local riemannian geometry,
$Spin^{\C}$ - structures and Seiberg - Witten equation
      \endtitle
  \author Nikolay Tyurin                    \endauthor

   \address BLTP JINR
   \endaddress
  \email
    tyurin\@tyurin.mccme.ru  ntyurin\@thsun1.jinr.ru
   \endemail

  \endtopmatter

\head Introduction
\endhead

The mathematical part of Seiberg - Witten theory is now a powerful
tool in the framework of smooth geometry of low dimensions. Using the
celebrated equations introduced by E. Witten in [1] one could establish
a number of new results and   reprove known ones. In the original paper
[1] the author first of all reproves Donaldson fundamental results
which ensure us that there are exist homeomorphic but not
diffeomorphic smooth compact 4 - dimensional manifolds. The key point is
that for any Kahler manifold one has nontriviality of the invariant
(which is called Seiberg - Witten invariant) while for almost all "topomodels"
the invariant is trivial (we'll give all the definitions and constructions
below). This fact itself destroys a number of physically - philosophical
speculations  about the structure of Universe;
usually they discussed what is the topological structure and what is the
manifold which corresponds (and represents) our world. But now if
say the famous Poincare smooth conjecture is not true (one can suppose
it because despite of the huge attempts to prove it we have not
yet any answer on this question) then the discussion would be
illuminated by new colors. The problem arose is not only what is
the manifold but additionally  which   smooth structure plays on. For example
usual $\R^4$ looks like a bad joke: Donaldson proved that it admits
continuum  set of different smooth structures. And it doesn't matter
what a metric one uses on this $\R^4$, euclidean or lorentzian.

But not only known results were proved using the new invariant. For
the Donaldson theory it was usual to compute the Donaldson invariants
in the framework of algebraic geometry. Almost all examples in
[23] come in this way. It was established by Taubes that
Seiberg - Witten invariants are related as well with symplectic
geometry --- in [3] he proves that for any symplectic manifold
with $b_2^+ > 1$ the associated canonical class is the basic one.
The proof was so simple and natural that  it was conjected that
the symplectic geometry is exactly the subject containing all
non trivial cases for the Seiberg - Witten invariants. This
suggestion was disproved almost immediately by Taubes himself
together with D. Kotschik   and J. Morgan in [5], but
the counterexample is a 4 - dimensional manifold with nontrivial
fundamental group. It means that one can correct the suggestion
just imposing the condition  of simply connectedness.

There are many courses now in the Seiberg - Witten theory
(we mention   in the References a number of these lectures
or surveys) so I hadn't in mind to present a textbook composed
of the materials which can be derived from the library.
I want just to add some geometrical facts and constructions
whose presence simplifying the reading of papers on the Seiberg -
Witten theory. The most attention will be payed to the notion
of $Spin^{\C}$ - structure and to the relation between complex
and spinor geometries. The second aim is to ask some natural questions
arise during the discussion which are from my point of view
meaningful mathematical problems themselfs. 

I would like to thank Korean Institute for Advanced Study (Seoul)
where the lectures were given, Joint Institute for Nuclear Research (Dubna)
where the writing was started and Max - Planck - Institute fur Matematik (Bonn) where I finished these remarks about this still unfinished story.

\head 1. The complex geometry and the spinor geometry
\endhead

Let $X$ be a smooth orientable compact riemannian manifold
of real dimension 4. Topological data, characterized $X$, are
its fundamental group $\pi_1(X)$, cohomology groups $H^i(X, \Z)$
whose ranks are usually called Betti numbers $b_i(X)$ of $X$ and
the intersection form $Q_X$ which is an integer quadratic form
on the lattice $H^2(X, \Z)$ whose definition requires the choice of
an orientation over $X$. Thus in simply connected case it remains
only the integer lattice and the intersection form encoding
the topological information about $X$. In the  gauge theories
(in the framework of differential geometry one usually combines the Donaldson
and the Seiberg - Witten theories using the term "gauge theories")  it's
quite normeous to consider only the case when the underlying manifold
$X$ has trivial fundamental group. Once again in this case one has
only lattice data.  The changing (or reversing)  of the orientation
multiplies the intersection form by $-1$. If an appropriate orientation
is fixed then one can reduce the intersection form $Q_X$
in enlarged space $H^2(X, \R)$ to a diagonal form and take the numbers
of positive and negative squares. These numbers are denoted as $b_2^+(X)$
and $b^-_2(X)$; the sum of these numbers is equal to the second Betti number
$b_2(X)$.

\subheading{Example} The simplest examples are $S^4$ and $\C \proj^2$ ---
4 - dimensional real sphere and complex projective plane. The last one
being considered as a 4 - dimensional real manifold is endowed with
the canonical orientation coming from the complex structure. Both the
4 - dimensional manifolds have trivial fundamental groups. The sphere
has trivial second Betti number so the lattice $H^2(S^4, \Z)$ and
consequently the intersection form are trivial. In the projective plane
one has a well known 2 - dimensional submanifold ---  projective line,
which is from the real geometry point of view just a 2 dimensional sphere.
We know from the standard algebraic geometry that every two projective
lines either coincide or have exactly one mutual point. This
intersection point is a "complex" point that is the index of
intersection is equal to 1. It follows from the positivity of
the complex  geometry. Thus the lattice $H^2(\C \proj^2, Z)$
is isomorphic to $Z[h]$ with the natural intersection form
$$
Q(k[h], l[h]) = k \cdot l.
$$
This intersection form is usually denoted as $1$.
In this case $b_2^+ = 1, b_2^- = 0$. Below we'll meet $\overline{\C
\proj}^2$ --- the projective plane with reverse orientation.
This manifold plays quite important role in the constructions.
It's clear that it has intersection form $-1$.

One has the tangent bundle $TX$ over our real based manifold which as a real
bundle has the following characteristic classes: the Shtiffel - Witney classes
$w_i$ belong to $H^i(X, \Z_2)$, the Pontryagin class $p_1 (X)$
belongs to $H^2(X, \Z)$ and the Euler class which gives us the euler
characteristic $\chi$ of our based manifold $X$. This euler characteristic
is related with the Betti numbers as follows
$$
\chi(X) = b_0(X) - b_1(X) + b_2(X) - b_3(X) + b_4(X) =
2(b_0(X) - b_1(X)) + b_2^+(X) + b_2^-(X).
$$
Here we consider the connected case so $b_0(X) = 1$. On the other hand
the Pontryagin class is related to the signature of the based manifold
$$
\si(X) = b_2^+ - b_2^- = - 3 \int_X p_1(X).
$$
Coming further let us suppose that our based manifold admits an almost
complex structure $I$ so is a nondegenerated operator
$$
I: TX \to TX
$$
with square
$$
I^2 = - id.
$$
Such an operator can exist even if $X$ doesn't admit a complex structure.
Here and below we establish that the existence of an almost complex structure
over $X$ is pure topological (or arithmetical) fact. Let us now examine the
case when such a structure exists. It means that the real bundle $TX$
endowed with $I$ is a complex rank 2 bundle $T^{1,0}X$ which has
as a complex rank 2 bundle the first and the second Chern classes
$c_1$ and $c_2$. The first Chern class is fixed modulo two by the
following condition
$$
c_1 = w_2(X) (mod 2)
$$ being an integer lifting of the second Shtiffel - Whitney class $w_2(X)$.
At the same time we can realize the Pontryagin class $p_1(X)$
by the definition as the second Chern class $c_2(T_{\C}X)$ of the
complexified tangent bundle $T_{\C}X$ (let us recall that  odd
Chern classes of any complexified bundle are trivial).  So
$$
T_{\C} X = T^{1,0}X \oplus T^{0,1}X = T^{1,0} X \oplus \overline{T^{1,
0}X}
$$
and it follows that
$$
c_2 (T_{\C} X) = 2 c_2  - c_1^2.
$$
The integration of the last formula over $X$ gives us the
following relationship
$$
c_1^2 = 2 \chi(X) + 3 \si(X),
$$
which is called Wu formula. The arithmetics comes with the statement
of Wu theorem which claims that the existence of an almost complex structure
over a 4 - dimensional manifold is equivalent to the existence of
solutions to the following pure arithmetical conditions:
$$
c_1 = w_2(X) (mod 2) \quad \quad {\text and } \quad \quad c_1^2 = 2 \chi
+ 3 \si.
$$
We've proved the theorem in one direction, using a brilliant result of
Hirzebruch, who proved that
$$
\int_X p_1(TX) = - 3 \si (X).
$$

\subheading{Example} Really it was unknown in 50th  is there any almost
complex structure over $S^4$ or not. Nowadays we can solve this problem
almost immediately following F. Hirzebruch. Let us compute the number
$$
\kappa = 2 \chi + 3 \si
$$
for 4 - dimensional sphere $S^4$. We have
$$
\kappa (S^4) = 2 (2 - 0 + 0) + 3 (0) = 4.
$$
At the same time $S^4$ has no nontrivial 2 - cohomological classes
with nontrivial squares. So the answer to the old problem is dictated
by this pure arithmetical fact.

At the same time for $\C \proj^2$ we have
$$
\kappa (\C \proj^2) = 2 (2 - 0 + 1) + 3 (1) = 9.
$$
Therefore one gets
$$
c_1 = \pm 3 [h]
$$
and this one is not an unexpected answer.

At the same time for the reverse oriented manifold
$\overline{\C \proj}^2$ the computation gives
$$
\kappa (\overline{\C \proj}^2) = 2 (2 - 0 + 1) + 3 (- 1) =
3,
$$
but every 2 - cohomological class in $H^2(\overline{\C \proj}^2, \Z)$
has nonpositive square.

Now we want to realize the numbers which come to the signature in
a different style. Namely let us fix a riemannian metric $g$ over $X$.
Together with the fixed orientation it gives us the Hodge star
operator:
$$
\aligned
*: \Om_X^i \to \Om_X^{4-i}; \\
<\be, * \al>_g  d \mu = \be \wedge \al \\
\endaligned
$$
for any $\be \in \Om_X^{4-i}$ where $d \mu$ is the volume form. This
operator could be exploited to deform the usual de Rham complex where
the role of ordinary $d$ plays the formal adjoint operator $d^*$ which is
a conjugation of $d$ by the Hodge star operator. A famous Hodge theorem
tells us that for every $i = 0, ..., 4$ one can decompose the space
of $i$ - forms into the following direct sum
$$
\Om^i = \Cal H^i \oplus d \Om^{i-1} \oplus d^* \Om^{i+1}
$$
and this decomposition is orthogonal with respect to the scalar product,
induced by the riemannian metric $g$. The first space $\Cal H^i$
consists of harmonic $i$ - form which satisfy
$$
\Cal H^i = \{ \rho \in \Om_X^i \quad | \quad d \rho = d^* \rho = 0\}.
$$
The space has finite dimension equals to $b_i$. So every cohomological class
is represented uniqually by a harmonic form.

Further, 2 - forms over 4 - dimensional manifold have a specialty in the
framework of the Hodge theory namely it's not hard to see that the Hodge
star operator being considered on $\Om_X^2$ is an involution
$$
*: \Om^2_X \to \Om^2_X, \quad \quad *^2 = id
$$
(it follows from the fact that the wedge product is symmetric
on $\Om^2_X$). Thus one can decompose the space of 2 - forms
into the direct sum of eigenspaces of the Hodge star operator.
Let $\Om^+$ be the eigenspace with eigenvalue $+1$ while $\Om^-$
corresponds to $-1$. It means that for every "positive" 2 -
form $\al \in \Om^+$ the wedge product is positive
$$
al \wedge \al = \vert \al \vert^2_g d \mu
$$
and for every "negative" $\be \in \Om^-$ the wedge square is negative
$$
\be \wedge \be = - \vert \be \vert^2_g d \mu.
$$
Moreover one takes the intersections
$$
\Cal H^{\pm} = \Cal H^i \cap \Om^{\pm}.
$$
The space $\Om^+$ is called the space of self dual forms while
$\Om^-$ is called the space of anti self dual forms. Since the decomposition
is pure local it induces the corresponding decomposition of
the bundle $\La^2 T^*X$ into the direct sum of two subbundles
$$
\La^2 T^*X = \La^+ \oplus \La^-.
$$
Let us recall that in dimension 4 all the decompositions and notions
are comformally invariant so they depend only on conformal class
of the riemannian metric rather then on riemannian metric itself.

As it was stated in a famous Hodge theorem every cohomological class
is represented by a harmonic form with respect to every riemannian metric.
But the case "2 in 4" is special, and one could ask whether or not a given
2 - cohomological class $[al] \in H^2(X, \Z)$ can be realized by a self dual
(or anti self dual) harmonic 2 - form. A priori we have only one type
of 2 - classes which definitely cann't be represented by self or anti self
dual forms namely   isotropical 2 - classes with zero squares. For
any 2 - cohomological class $\eta \in H^2(X, \Z)$ with positive square
one has an open question in smooth 4 - geometry: is there exist such
riemannian metric (more rigouresly - such conformal class of riemannian
metrics) that the harmonic with respect to this metric form $\al_{\eta}$,
representing the class, is self dual. The same question takes place
for any 2 - cohomological class with negative square. These two questions
are closely related since the reversing of the orientation interchanges
the notions of self duality and anti self duality. But we want to formulate
an important unsolved problem in 4 - dimensional differential geometry
in more general and geometric form. To do this we need to introduce
the so called real period map.

Let $X$ be a smooth compact orientable riemannian 4 - manifold with signature
$$
\si = b_2^+ - b_2^-.
$$
The standard slant product (which gives us the intersection form $Q_X$)
extended to the space $H^2(X, \R)$ gives us a bilinear pairing
on this real vector space. We will call a $b_2^+$ - dimensional
subspace $V \subset H^2(X, \R)$ positive if the restriction of this
bilinear pairing is strictly positive on $V$. At the same time
one has the notion of negative subspace of dimension $b_2^-$
but for every positive $V^+$ one has unique negative $V^-$
taking the  orthogonal compliment of the first one. Now let
us take the positive Grassmann manifold $Gr^+(H^2(X, \R)$
consists of all positive $b_2^+$ subspaces. It is
canonically isomorphic to the negative Grassmann $Gr^-(H^2(X, \R)$
so one can remove the choice of orientation identifying these
Grassmann manifolds. If we choose a riemannian metric (or
a conformal class) over $X$ we get two points $p^{\pm} \in Gr^{\pm}
(H^2(X, \R))$ which correspond to the subspaces of self dual and anti self
dual harmonic 2 - forms respectively. Again if we forget about the
choice of  orientation then we get identifying the Grassmann manifolds and
the points a map
$$
P: \Cal R \to Gr^+(H^2(X, \Z)),
$$
where $\Cal R$ is the space of all riemannian metric  the given
smooth structure on $X$ admits
which is called real period map. It is a real analogy
of the celebrated period map in algebraic geometry.

S. Donaldson (see [23]) proved that this map is locally surjective.
It means that for every point $p \in Gr^+(H^2(X, \R))$ realized by
say self dual harmonic subspace in $\Om^2_X$ there exists a small
neighborhood consists of points which can be realized as well.
But generally one has the following important:

\proclaim{Problem} Is the real period map  globally surjective?
\endproclaim

It would be strange if it is not so but nevertheless nobody knows
now is it true or not. I have to say here that this problem
is in some sense an analogy of the famous Hodge conjecture in
algebraic geometry. These two problems are quite parallel:
in the algebraic geometrical setup the condition that a
cohomological class has (l, l) - type is of the same level as in
the differential geometrical setup the condition that a cohomological class
lies on a positive subspace. So if we thought of these two conditions
as analogous conditions (cohomological conditions) then these two
problems are  about realizations of these classes by real objects:
some algebraic subvariety and some harmonic self dual (or anti self
dual) form.

Thus we've discussed what a riemannian metric chosen over $X$ brings to
the picture.  But it is only one part; let us turn to the second.
We call this part "twistor and spinor geometry". As we'll see it
is extremely close to complex geometry but usually isn't
considered in the framework of the last one.

Let one take a riemannian metric $g$ over our based smooth manifold
$X$ and fix an appropriate orientation. Then one gets two projective
bundles $\proj^{\pm} \to X$ which is called twistor bundles.
We would like to give below two description of these bundles
having in mind to relate closely spinor and complex geometries.
The first description is quite classical: over each point
$x \in X$ one takes the set of all complex structures which are
compatible with the given riemannian metric $g$. Identifying
the fiber $T_xX$ of the tangent bundle with standard $\R^4$
and reduce the riemannian metric to the standard form one gets
the standard description of the set in terms if operators $I, J, K$
(the standard quaternionic generators) so any compatible
complex operator can be represented as
$$
A = \al I + \be J + \ga K, \quad \quad \al, \be, \ga \in \R
$$
such that
$$
\al^2 + \be^2 + \ga^2 = 1.
$$
Thus one gets a 2 - dimensional sphere  which is topologically
isomorphic to  projective line. Now one can change the given orientation
and consider the complex structures which are compatible with $g$
and with the reverse orientation. This gives us again projective line.
Let us distinguish these two lines assigning $+$ to the first one and
$-$ to the second. Generalizing this local picture over all $X$ one
gets two projective bundles
$$
\proj^{\pm} \to X,
$$
such that each fiber consists of complex operators compatible with
our metric and given (or reverse) orientation. This definition is
the classical one. The total space of the projective bundle
$\proj^- \to X$ is called "twistor space" (the choice
of the sign is conventional; for this twistor space
one has the wonderful relationship between instantons
over $X$ and stable holomorphic vector bundles over
the twistor space, see [12]). It is
a real 6 - dimensional manifold fibered over $X$ on 2 - dimensional
sphere. Moreover this 6 - dimensional real manifold is canonically
endowed by an almost complex structure.
Since we want to formulate a relative problem in a future
let us recall how one gets this almost complex structure.

Every fiber has its own complex structure being a complex
projective manifold (we will show in a moment why the 2 -
sphere has the fixed structure). Our riemannian metric
$g$ defines the Levi - Civita connection which induces a
special connection on the projective bundle $\proj^-$.
If we denote as $Y$ the total space of the projective bundle
then the tangent space $T_y Y$ in arbitrary point $y \in Y$
is decomposed due to the induced connection into two
direct summands
$$
T_yY = V \oplus W
\tag *
$$
where $V$ is the tangent space to the fiber in this point and
$W$ is isomorphic to $T_xX$ where $x$ is the image of $y$
with respect to the canonical projection
$$
\pi: Y \to X.
$$
Note that the first summand $V$ has a complex operator coming from
the complex structure on the fiber. On the second summand $W$ one takes
the complex operator which corresponds to the point $y$ of the fiber
(all points in the fibers are complex operators on $T_xX$). Thus
the direct sum of the operators gives us a complex operator on
$T_yY$. Globalizing this picture one gets an almost complex structure
over all $Y$. This is really twisted structure since it isn't
constant along the fibers but is twisted along these.
It's a well known fact (see f.e. [12]) that the integrability condition
for this almost complex structure is satisfied if and only if the started
riemannian metric is anti self dual.

\subheading{Example} There are two mostly used examples of integrable
twistor spaces: the first one is for $S^4$ with the standard conforamally
flat riemannian metric (the corresponding twistor space
is $\C \proj^3$) and the second is for $\C \proj^2$ with the standard
Kahler metric (here one gets the flag variety, see [12]).

Now let us discuss the second description of the twistor bundles.
For this again consider the riemannian case locally. Let $x \in X$
is a point of $X$. The tangent space $T_xX$ is a real space endowed with
symmetric positive bilinear form induce by fixed riemannian metric $g$.
Let complexify the space $T_xX$ and then projectivize it. During the
complexification we extend $g_x$ on $T^{\C}_xX$ as a complex bilinear
form. So we get on the complex space $T^{\C}_xX$ not
a hermitian form but a complex metric (in the style of [14]).
This complex bilinear form $g_x^{\C}$ gives us after the
projectivization a nondegenerated quadric
$$
Q_x \subset \C \proj^3 = \proj (T^{\C}_xX),
$$
consists of isotropical vectors of $g_x^{\C}$. Moreover one has
a standard real structure $\Theta$ on this $\C \proj^3$   inherited
during the complexification (since our space was real before). It's not
hard to see that

- our quadric $Q_x$ is a real quadric since our metric $g$ was real
before complexification

and

- it hasn't real points at all since our metric $g$ being riemannian
doesn't admit any real isotropic vectors.

Thus locally riemannian geometry in dimension 4 gives
at each point the following  data: $\C \proj^3$ with a standard
real structure and a real quadric $Q \subset \C \proj^3$ without
real points.

It's well known fact (see f.e. [13]) that every nondegenerated quadric
in 3 - dimensional projective space is the direct product
of two projective lines
$$
Q = \proj^+ \times \proj^-.
$$
One can understand it more geometrically: there are two family of
projective lines lie on the quadric and $\proj^+$ and $\proj^-$
parameterize the families.
So points of $\proj^+$ one can understand as projective lines
lie on $Q$ and belong to the first family while $\proj^-$
consists of projective lines from the second one. One
can distinguish these two family due to the fixed orientation.
As we will see the first family just corresponds to
complex operators compatible with the given orientation while
the second one is associated with reverse orientation.
To see this more explicitly let us study what is a compatible
complex structure in this local picture.

Let $I$ be a local almost complex structure, compatible with
the given riemannian metric and orientation. When we complexify
the tangent space $T_xX$ we can decompose the complexified space
into the direct sum of two part with respect to by - types:
$$
T_x^{\C} X = T^{1,0}_x \oplus T^{0,1}_x,
$$
where the first summand is the eigenspace of $I$ with eigenvalue
$\imath$ and the second one is the eigenspace with eigenvalue
$- \imath$ as usual. One can see that this decomposition defines
our operator $I$ uniquely. After the projectivization this gives us
two projective lines in our $\C \proj^2$ namely
$$
\proj(T^{1,0}_x) , \proj(T^{0,1}_x) \subset \C \proj^3.
$$
Denote these lines as $l_I$ and $\bar l_I$. We know which one is
the first and which one is the conjugation of the first one
by our standard real structure $\Theta$:
$$
\Theta(l_I) = \bar l_I.
$$
It's clear that these two projective lines have no mutual points.
Therefore local picture for complex geometry in dimension 4 is
as follows:  the projective space $\C \proj^3$ together with
a standard real structure and a pair of conjugated projective lines
without real points.

The next step is to combine these two local picture: what is a riemannian
metric with a compatible almost complex structure? The compatibility
condition for $g$ and $I$ is quite simple: riemannian metric
$g$ and complex structure $I$ are compatible if and only if
at each point of $X$ the projective lines $l_I, \bar l_I$ lie
on the quadric $Q$. But since the intersection of
$l_I$ and $\bar l_I$ is trivial these two projective lines
belong to one family. If $I$ is compatible with given orientation
then $l_I, \bar l_I \in \proj^+$. If it induces reverse orientation
then $l_I, \bar l_I \in \proj^-$. Therefore the projective line
$\proj^+$ parameterizes complex structures compatible with given
orientation and $\proj^-$ respects for almost complex structures
induce reverse orientation. Thus we get here two 2 - spheres
which were called above "twistor bundles". But now we get
much more important information about these
twistor bundles, namely from the local picture in $\C \proj^3$
we derive that these $\proj^{\pm}$ are endowed with their
own complex structures. Moreover, our projective bundles are
$PU(2)$ - bundles. They automatically carry not only own
complex structures but as well every fiber has a kahler structure!
To see this explicitly we have to add to the picture the action
of our standard real structure $\Theta$. Since our quadric is real
this means that one has two possibilities for the families of projective
lines on it: either $\Theta$ maps one family to the other or
$\Theta$ maps each family to itself. The first case is forbidden because
of riemannian property of our metric. If $\Theta$ maps $\proj^+$ to
$\proj^-$ then it has to be real points on $Q$! Really, if one takes
a line $l$ say from  $\proj^+$ and if $\Theta(l)$ belongs to $\proj^-$
then it would be the intersection point $p = l \cap \Theta(l)$
and it should be real. It implies that
$$
\Theta(\proj^{\pm}) = \proj^{\pm}.
$$
Further, the action of $\Theta$ on $\proj^{\pm}$, preserving the lines,
induces a pair of quaternionic real structures $\theta^{\pm}$
on the lines. Let us recall that a quaternionic real structure
as an anti holomorphic involution on   $\C \proj^m$ where $m$ is odd
which has the following form in  appropriate homogeneous coordinates
$$
\theta (z_0: z_1: ...: z_{m-1}: z_m) =
(- \bar z_1: \bar z_0: - \bar z_m: \bar z_{m-1})
$$
(while a standard real structure exists on any projective space and
has the form
$$
\Theta(z_0: ...: z_m) = (\bar z_0: ... : \bar z_m)).
$$
Of course the points $l_I$ and $\bar l_I$ are conjugated by $\theta^+$.

It is an additional ingredient which we didn't use yet: we've considered
the bilinear extension of the riemannian metric $g$ to the complexified
space $T^{\C}_xX$ but as well one has the sesqueliniar extension
makes the complex space $T^{\C}_xX$ into a hermitian space and consequently
makes the projective space $\C \proj^3$ into a kahler manifold
endowed with a kahler metric with a kahler form. Of course
this kahler metric by the construction is invariant under
the action of our standard real structure $\Theta$. Therefore it
splits with respect to the decomposition
$$
Q = \proj^+ \times \proj^-
$$
and induces two kahler metrics on $\proj^{\pm}$. We have to emphasize
here that these kahler structures depend on riemannian metric rather then
the previously described objects depend only on its conformal class.
It is the striking difference between  conformal and riemannian geometries
over twistor spaces. To illustrate the difference we would like to propose
a problem which goes quite parallel with the classical solved problem
mentioned above about integrabily of the twistor almost complex structure
over a twistor space. Now let us emphasize again that we fix a riemannian
metric on our based manifold $X$. Then on the twistor space $Y \to X$
we have a hermitian triple $(G, I, \Om)$ consists of a riemannian metric
$G$, a compatible almost complex structure $I$ and the corresponding
almost kahler form $\Om$. The second element is the same as in the classical
twistor construction and has been described above. To construct the triple
one can either describe the first element or the third one because
in a hermitian triple every two elements uniquely define the third
element. Using again the decomposition (*) induced by the Levi -
Civita connection one takes at the point $y \in Y$ the following
direct sum of the kahler form on the fiber and the "twistor" 2 - form
on the horizontal component which corresponds to the given riemannian metric
$g$ and the complex structure operator represented by the point of the
fiber. Totally we get a nondegenerated globally defined 2 - form $\Om$
on $Y$ which has type (1,1) with respect to the twistor almost
complex structure $I$ and it is a positive form. This means that
it together with $I$ defines a riemannian metric $G$ over the
twistor space $Y$ and thus the hermitian triple $(G, I, \Om)$ is constructed.

Thus one could say that the difference between conformal and riemannian
geometries can be formulated on the twistor space level as the difference
between almost complex and almost kahler geometries. It implies
the following very natural and important question.

\proclaim{Problem} When the induced by a riemannian metric $g$
almost kahler form $\Om$ on the twistor space $Y$ is closed
so is a symplectic form?
\endproclaim

The solution to this problem would open a way to get the twistor
construction for the Seiberg - Witten theory. Let us recall that
for instantons the twistor construction gives us an important
relationship between anti self dual connections and holomorphic
vector bundles over the twistor space. So the paradigm
for the Donaldson theory is the geometry of stable
holomorphic vector bundles. At the same time
there is now  a very popular idea which tries to
explain the mirror duality  comparing
holomorphic and symplectic geometries. If one has two
mirror manifolds then it should be a map from the holomrphic
objects of the first one to a symplectic objects of the second
manifold and vice versa. On the other hand there is a
miracle duality between the Donaldson theory and the Seiberg -
Witten theory which was explained physically by E. Witten and
mathematically by A. Tyurin and V. Pidstrigach. But the explanations
give us only numerical duality so as the values of the invariants
coincide. But if one believes that more strong then pure numerical
duality takes place it should be a construction relates
the theories on the twistor level, and this duality on
the twistor level would be reflected by a duality
between holomorphic and symplectic setups. It is
just a suggestion now but anyway the problem formulated above
is quite interesting and meaningful itself.

Now the time of spinor geometry comes. Since our based manifold is orientable
then according to an old classical result of Hirzebruch and Hopf
every projective bundle can be lifted to a vector bundle.
Such a lifting depends on the choice of the first Chern class.
Namely if we choose a class $c \in H^2(X, \Z)$ such that
$$
c = w_2(X) (mod 2)
$$
where $w_2(X)$ is the second Shtiffel - Whitney class $w_2(X) \in H^2(X, \Z_2)$
of $X$ then we automatically get the corresponding lifts of $\proj^{\pm}$
to a pair of $U(2)$ - bundles $W^{\pm}$ with the first Chern classes
$$
c_1(W^+) = c_1(W^-) = c.
$$
This choice is called a choice of $Spin^{\C}$ - structure on $X$.
For every such lifting the pair $W^{\pm}$ are endowed automatically
by hermitian structures induced by our riemannian metric $g$
and satisfy the following fundamental for the theory properties:

1) the homomorphism bundle $Hom (W^-, W^+)$ is isomorphic
to $T^{\C}X$, moreover, the "real" part $Hom_J (W^-, W^+)$
consists of the maps which preserve both hermitian structures
is isomorphic to $TX$;

2) the adjoint bundles $ad W^{\pm}$ are isomorphic to $\La^{\pm}$
modulo diagonal parts
where $\La^{\pm}$ are the direct summands in the decomposition
of $\La^2 T^*X$ induced by the Hodge star operator.

Let us recall that for a given hermitian vector bundle one can
associate the corresponding principle bundle and then represent the
last one using the adjoint representation. The resulting
adjoint bundle for $W^{\pm}$ is locally defined by skew hermitian
endomorphisms. If $\phi$ is a spinor field so is a section of
say $\W^+$ then one has using the hermitian structure
a endomorphism
$$
a_{\phi} = \phi \otimes \bar \phi \in End W^+,
$$
which preserves the hermitian structure. If our spinor bundle
$W^+$ moreover is $SU(2)$ - bundle so the determinant line bundle
$det W^+$ is trivial then $ad W^+$ is exactly isomorphic  to
$\La^+$ since one has a local isomorphism of $SU(2)$ and $SO(3)$.
This case is special; this means that our chosen $Spin^{\C}$ - structure
$c \in H^2(X, \Z)$ is trivial and it's possible if and only if
the second Shtiffel - Whitney class $w_2(X)$ is trivial. One calls
the manifolds with trivial $w_2(X)$ spin manifolds. The distinguished
lifting of the twistor bundles to $SU(2)$ - bundles is called
the choice of $Spin$ - structure. And repeating again in  this
spin - case one has exact isomorphisms
$$
ad W^{\pm} \equiv \La^{\pm}.
$$
We would like to generalize these original relationships
to the case when $W^{\pm}$ have non trivial determinant.
For this in the framework of the Seiberg - Witten theory it's
sufficient to construct the following quadratic maps
$$
q: W^{\pm} \to \La^{\pm}
$$
based  again on complex geometry. Namely as we have seen
for our given metric  $g$ one can define compatible
local almost complex structures in terms of spinor geometry.
If $\phi \in \Ga (W^+)$ is a section which doesn't
vanish in a neighborhood $O(x)$ of a chosen point $x \in X$
then it defines a compatible almost complex structure
$I_{\phi}$. Namely the projectivization of the section
gives us a local section of the twistor bundle $\proj^+$
and this is clearly as we have seen above an almost complex structure
defined in the neighborhood $O(x)$. Thus one has the corresponding
self dual 2 - form $\om$, reconstructed from $g$ and $I_{\phi}$.
But this is a feature of twistor geometry; for spinor geometry
all the story is lifted to the vector level therefore
one has a correspondence between two multiplications:
every spinor field can be multiplied by a local complex function
while every self dual form can be multiplied by a local real function.
Then if we realize the corresponding section of the twistor bundle
as a normalized spinor field with a fixed phase then the multiplication
by a complex function $z$ of the spinor field induces
the multiplication by $z \bar z$ of the distinguished self dual 2 - form.
Since $det W^+$ is not trivial then one couldn't choose an apporpirate phase
fixing simultaneously over $X$; nevertheless if one takes the bundle
$ad W^+$ as a subbundle of $End W^+$ consists of skew hermitian endomorphisms
then going further one can take inside of $End W^+$ the projection
of $ad W^+$ to the subspace of traceless endomorphisms. This gives us
the possibility to compare spinor fields with self dual forms;
reversing the orientation one gets the same story for $ad W^-$ and $\La^-$.

At the same time we can continue our discussion of the Wu theorem which we've
proved in one direction. Now to complete the proof let us just mention here
that the question of the existence of a compatible almost complex structure
for a given riemannian metric $g$ is pure topological. Namely if there exists
a lifting  of the twistor bundle $\proj^+$ to such spinor bundle $W^+$
that it admits a nonvanishing smooth section $\phi_0$ then there exists
 globally defined compatible almost complex structure. And such a section
exists if and only if the second Chern class $c_2(W^+)$ is trivial.
Therefore if a numerical condition holds for a $Spin^{\C}$ - structure
$c = c_1(W^{\pm}) \in H^2(X, \Z)$ then there exists a nonvanishing spinor
field and consequently a globally defined almost complex structure.

Using these arguments together with some additional remarks one can
establish that if one chooses an appropriate class $c \in H^2(X, \Z)$
to lift $\proj^{\pm}$ to vector bundles $W^{\pm}$ then determinently
$$
\aligned
c_2(W^+) = \frac{1}{4}(c^2 - 2 \chi - 3 \si), \\
c_2(W^-) = c_2(W^+) + \chi, \\
\endaligned
$$
where as usual $\chi$ is the euler characteristic of $X$ and
$\si$ is the signature. Thus one should see that topologically
$Spin^{\C}$ - structure $c \in H^2(X, \Z)$ is a generalization
of the notion of canonical class.

To finish this lecture let us illustrate the properties of spinor
bundles mentioned above on the projective level.

First of all, we have $TX = Hom_J (W^-, W^+)$. Returning to
our $\C \proj^3$ endowed with a standard real structure $\Theta$
and containing a nondegenerated  real quadric $Q$ without real points
we have to represent each point of $\C \proj^3$ as a map from $\proj^-$
to $\proj^+$ where $\proj^{\pm}$ as usual parameterize two families
of projective lines on $Q$. The map $\psi_p$ is defined as follows:
for every projective line $l \in \proj^-$ you take the projective plane
$\pi(p, l)$ which is spanned by the point $p$ and the line $l$.
This plane intersects the quadric $Q$ in a plane curve of dimension 2
as it is predicted by Bezout theorem. But this intersection curve has
 a special type since $l$ itself lies on $Q$. This means that
the intersection consists of exactly two projective lines:
line $l \in \proj^-$ and a projective line $\psi_p(l)$
which intersects $l$ and lies on $Q$ therefore it contains
in $\proj^+$. Thus we defined the map
$$
\psi_p: \proj^- \to \proj^+
$$
for each point $p \in \C \proj^3$ say away from the quadric $Q$. Really we
want to construct these maps only for real points of $\C \proj^3$ which
correspond to $TX$ rather then to $T^{\C}X$. It means that we can leave
the question about the maps for the points of $Q$ turning to
the reality question for $\psi_p$ if $p$ is real. But this question is
almost trivial: if $p$ is a real point of $\C \proj^3$ then it remains
unchanged during the action of our standard real structure $\Theta$.
But as we mentioned above the action induces two quaternionic real structures
on $\proj^{\pm}$. Then we get applying $\Theta$ to $\psi_p$ the following
commutation relation
$$
\psi_p (\theta^- (l)) = \theta^+ (\psi_p (l))
$$
for real $p$ and this gives the reality condition which declare
that $\psi_p$ for real $p$ preserves both the quaternionic structure.

To compare $ad W^{\pm}$ and $\La^{\pm}$ we can consider another projective
space coming from local consideration. Namely let
one takes the fiber of $\La^2 T^*X$ over the point $x \in X$, complexify
it and then projectivize. Since $\La^2 T^*_xX$ has real dimension
6 then one gets $\C \proj^5$. In this projective space we have
the Grassmann quadric $Gr \subset \C \proj^5$ which paramertizes
all projective lines in our $\C \proj^3$. Really for every complex 4 -
dimensional space $V$  the set of 2 - dimensional subspaces can be
identified with a special class of 2 - dimensional form up to constant  from
$\La^2 V^*$ which have rank 2. This condition is equivalent to
the equality
$$
\om \wedge \om = 0,
$$
and the last one is a quadratic equation whose solutions form Plucker quadric
$Gr$ in $\C \proj^5$. Thus the quadric $Gr$ represents all projective lines
in our $\C \proj^3$. Therefore one has two curves $l^{\pm} \subset Gr$
correspond to families $\proj^{\pm}$ defined by our riemannian metric $g$.
On the other hand the Hodge decomposition
$$
\La^2 T^*X = \La^+ \oplus \La^-
$$
gives us two projective planes
$$
\pi^+, \pi^- \subset \C \proj^5,
$$
and it's not hard to see that
$$
l^{\pm} = Gr \cap \pi^{\pm}.
$$
From this one sees that $l^{\pm}$ are two plane conics with involutions
without fixed points induced by quaternionic real structures $\theta^{\pm}$.
Now for every point $\phi \in l^+$ we can construct a class of self dual
2 - forms as follows. For the point $\phi$ one takes the conjugated point
$\bar \phi \in l^+$ and constructs two tangent to the conic $l^+$
projective lines $r_{\phi}, r_{\bar \phi} \subset \pi^+$ in the points.
There are no bi - tangent  lines to any conic hence for any $\phi$
the lines $r_{\phi}, r_{\bar \phi}$ are distinct and define uniquely
the intersection point
$$
\om = r_{\phi} \cap r_{\bar \phi}.
$$
This point $\om$ represents a class of 2 -forms which are defined
by the point up to a scalar multiple. So on the porjectiviztion
level one gets a quadratic map from spinor fields to self dual
2 - forms compatible with hermitian structures.

To finish we add that the same one takes place for spinor fields
from $W^-$ and anti self dual forms.

\head  2. The Seiberg - Witten equations
\endhead

The background idea of the previous lecture was that spinor and complex
geometries are very closed and even coincide in same points.
A background idea of the present one is  that the nature of the
Seiberg - Witten invariants is complex to.

Again let $X$ be a real 4 - dimensional compact orientable smooth
riemannian manifold and $g$ is a fixed riemannian metric on it.
Additionally we fix an orientation and choose a $Spin^{C}$ -
structure $c \in H^2(X, \Z)$ which is a lifting of projective
twistor bundles $\proj^{\pm}$ to a pair of spinor  $U(2)$
- bundles $W^{\pm}$ with the same determinant line bundle
$det W^{\pm}$ and hence with the same first Chern class
$c_1(W^{\pm}) = c$. These spinor bundles $W^{\pm}$ satisfy
as usual the properties listed in the previous lecture.

The configuration space for the system of equations which were
called the Seiberg - Witten equations is the following:
it is the direct product of $\sA_h(det W^+)$ --- the space of
hermitian connections on the determinant line bundle ---
and $\Ga(W^+)$ --- the space of all spinor fields that is
smooth sections of $W^+$. According to the physical tradition
we call the first component bosonic while the second fermionic
parts of the configuration space. For a point
$$
(a, \phi) \in \sA_h(det W^+) \times \Ga(W^+)
$$
of the configuration space one defines two equations. First of all
the Levi - Civita connection defined by our riemannian metric
induces a hermitian connection on $ad W^+$; at the same time one has
 the following usual decomposition
$$
\sA_h(W^+) = \sA_h (det W^+) \times \sA_h (ad W^+)
$$
hence for any hermitian connection $a \in \sA_h(det W^+)$
one automatically has the corresponding hermitian
connection on  whole bundle $W^+$ with covariant derivative
$\na_a$:
$$
\na_a: \Ga (W^+) \to \Ga (W^+ \otimes T^*X).
$$
But for $T^*X$ we have the conjugated identification
$$
T^* X = Hom (W^+, W^-)
$$
coming from the first property of $Spin^{\C}$ - structure.
Then one gets the following combination
$$
\aligned
\na_a: \Ga(W^+) \to \Ga(W^+ \otimes T^* X) = \Ga(W^+ \otimes Hom (W^+,
W^-)) = \\
\Ga (W^+ \otimes (W^+)^* \otimes W^+) \to \Ga (W^-) \\
\endaligned
$$
where at the last step we use the natural cancellation which is called
in this case Clifford multiplication. This combination gives us
an operator
$$
D_a: \Ga(W^+) \to \Ga(W^-)
$$
which is called Dirac operator coupled by the abelian connection $a$.
In a simple case when the based manifold is a spin manifold so admits
a canonical lifting of  twistor bundles  for a riemannian metric $g$
to spinor bundles with trivial determinants one has a preffered
Dirac operator on the spinor bundle $W^+_0$. Namely on the trivial
determinant bundle one has the trivial connection which corresponds
to ordinary differential and hence one takes the Dirac operator defined
by this trivial connection. This Dirac operator is a distinguished
operator but in our general $Spin^{\C}$ - case there is no any
preference to choose an operator. Originally Dirac got it
taking a "square root" from the Laplacian defined by
riemannian metric. The Laplace operator is a real operator
but complex geometry comes when one tries to define the square
root. As we will see the Dirac operator is exactly the square root
in the flat case; in a curved case with $Spin^{\C}$ - structure
one has to add the curvature of the abelian coupling connection
$a$ and the scalar curvature of the riemannian metric.

Thus the first equation of the Seiberg - Witten system reads as
$$
D_a(\phi) = 0
$$
and can be understood as the harmonicity condition for
spinor field with respect to coupled Dirac operator.

The second equation of the system is based on the second property of
$Spin^{\C}$ - structure. Following it one can compare the self dual part
of the curvature tensor of our abelian connection $a$ and a skew symmetric
traceless endomorphism of the spinor bundle $W^+$ defined by  spinor field
$\phi$ which we denote as $(\phi \otimes \bar \phi)_0$ having in mind
the projective story of the previous lecture. So the last equation reads
as
$$
F_a^+ = - (\phi \otimes \bar \phi)_0
$$
and now we can collect
these two together getting
$$
\aligned
D_a(\phi) = 0 \\
F_a^+ = - (\phi \otimes \bar \phi)_0
\\
\endaligned
$$
which is the celebrated Seiberg - Witten system of equations.
These equations are called the equations of abelian monopoles. Although
one defines a hermitian connection of whole $W^+$  the bosonic
variable is just abelian. On the other hand the gauge group of
the theory is abelian. It means as usual that one has a natural
group action on the configuration space such that the equations are
invariant under this action. In the case we are discussing
now such a group comes if we take
$$
\sG = Aut_h (det W^+)
$$
the group of fiberwise automorphisms of the determinant line bundle
which preserve its hermitian structure. So over each point on has
$U(1)$ as the automorphism group thus locally the gauge group
$\sG$ is modeled by maps
$$
\eta: N_x \to U(1)
$$
where $N_x$ is a small neighborhood of a point $x \in X$. On the first
component of the configuration space it acts in the usual way
while in the second case it acts only on the central part
(so locally if $\eta$ is the multiplication by a complex function
$e^{\psi}$ then $\eta (\psi) = e^{\psi} \phi$ where $\phi \in \Ga(W^+)$).
Thus one defines two quotient spaces: the first one comes with
factorization of the configuration space by the gauge group action
$$
\Cal C = \sA_h(det W^+) \times \Ga(W^+)/\sG;
$$
at the same time it's convenient to consider the subspace of irreducible
pairs $(a, \phi)$ which  have maximal stabilizer in the gauge group.
In our theory it just means that $(a, \phi)$ is reducible if and only
if the spinor component  $\phi$  vanishes identically. We define
the smooth part of the quotient space
$$
\Cal C^* = (\sA_h(det W^+) \times \Ga(W^+))^*/\sG
$$
where $(\sA_h(det W^+) \times \Ga(W^+))^*$ consists of irreducible
pairs.

The key point is that the Seiberg - Witten equations are invariant under
the gauge group action.  Therefore this invariance ensures that
one can consider the moduli space of solutions. "Moduli" means
that it is the set of solutions {\it modulo} gauge transformations.
We denote it as $\sM(g, c)$ and call it the moduli space of abelian
monopoles. The notation reflects the dependence on riemannian metric
$g$ and $Spin^{\C}$ - structure $c$.

The first question we would like to discuss is the problem
of reducible solutions. We will see that such a solution exists if
and only if the fixed riemannian metric $g$ satisfies a special property.

Namely let $(a, \phi)$ be a reducible solution
$$
(a, \phi) = (a, 0).
$$
Then the first equation of the systems is satisfied automatically while the
second gives us the condition
$$
F_a^+ = 0
$$
which is the equation of abelian instanton. Such an instanton
exists if and only if the fixed riemannian metric is of a special type:
the harmonic form representing the first Chern class $c_1 (det W^+) =
c$  is anti self dual. Really if $F^+_a$ is trivial then
$$
F_a = F^-_a,
$$
but any abelian curvature is a pure imaginary closed 2 - form
thus $d F^-_a = 0 = d^* F^-_a$ since every closed (anti) self dual 2 - form
is automatically co - closed and therefore harmonic. It means that
$$
\frac{i}{2 \pi} F_a
$$
is a real harmonic anti self dual 2 - form, but according to
the Chern - Weyl theory is represents the first Chern class
of the determinant line bundle $det W^+$ since our abelian connection
lives on this bundle. But this property is very special on riemannian metrics
(see discussion on real period maps from the previous lecture). For example
if $c^2 > 0$ there is no such riemannian metric at all!Really for the
existence of such a representation by an appropriate riemannian metric
it's necessary  that the square of the class is non positive
because there are no anti self dual forms with positive squares at all.
Moreover if $c^2 = 0$ and $c$ is non trivial element in the lattice
$H^2 (X, \Z)$ then there are no anti self dual representations of
this class because any anti self dual class with zero square has
to be trivial. Further even if there are exist riemannian metrics
admitting for $c$ anti self dual representations they form  sufficiently
"thin" subset in the space of riemannian metrics. Whole the space
is isomorphic to an infinite dimensional ball and the metrics admitting
such representations form a subset of codimension $b^+_2(X)$ (it follows
from the local surjectivity  theorem for the real period map proved
by S. Donaldson). From these notes one can derive that:

if $b_2^+(X)$ is bigger then zero then for a generic riemannian metric there
are no reducible solutions at all;

if $b^+_2(X) > 1$ then for every pair of generic metrics $g_0, g_1$
there exists a smooth path $\ga$ in the space of riemannian metrics
with ends at $g_0$ and $g_1$ such that every riemannian metric which is
parameterized by the path the harmonic 2 - form representing the class
$c$ is not anti self dual. So in these cases the moduli spaces are smooth
orbifolds.

The next question is what are the dimensions of moduli spaces in
these regular cases? To establish the numerical characterization
one has to study the linearization of the system. The first equation
$$
D_a(\phi) = 0
$$
is linear itself; from this we take the index of the coupled Dirac
operator. It is not very hard to compute this number applying
the well known Atiyah  - Hirzebruch index theorem which says that
every $\C$ - linear elliptic operator over a real compact smooth
manifold has the index which equals to a combination of
Chern classes of the "source" and "target" vector bundles and
$A$ - genus of the based manifold. Every algebraic geometer knows
this index theorem under the name of Riemann - Roch formula ---
and ommiting the calculations we claim that in our case
$$
ind D_a = \frac{1}{4}(c^2 - \si(X))
$$
(here we multiply the complex index by 2 to get the real one).
The second equation can be reduced to the infinitesimal level up
to gauge transformations with the following operator
$$
\de_a = d_a^* \oplus d_a^+: \Om^1_X(i \R) \to \Om^0_X(i \R) \oplus \Om^+_X(i \R)
$$
which is very well known in the Donaldson theory (see [23]). Its
index coincides with the index of the basic operator
$$
\de = d^* \oplus d^+: \Om^1_X \to \Om_X^0 \oplus \Om^+_X
$$
which defines "a half" of the standard de Rham complex. It gives
$$
ind \de_a = ind \de = - (1 - b_1(X) + b_2^+(X)) = - \half(\chi(X)
+ \si(X)).
$$
Thus the total index is
$$
ind SW = \frac{1}{4}(c^2 - \si) - \half (\chi + \si) =
\frac{1}{4}(c^2 - 2 \chi - 3 \si),
$$
therefore the virtual dimension of the moduli space $\sM(g, c)$
equals to
$$
v.dim \sM(g, c) = \frac{1}{4}(c^2 - 2 \chi - 3 \si).
$$
This number is very crucial for the theory at all: we have seen in
the previous lecture that this number is the second Chern
number of the positive spinor bundle
$$
v.dim \sM(g, c) = c_2 (W^+);
$$
if this number vanishes then $c$ is the canonical (or
anti canonical) class of some compatible almost complex structure.
For a generic metric $g$ if $b_2^+(X) > 0$ the moduli space is
a smooth real manifold of this dimension.

Until the time the story comes quite parallel to the consideration
in the framework of the Donaldson theory. As in the Donaldson case
we take a configuration space, consider a gauge invariant equation
and then construct the moduli space of solutions. The same question
of reducible solutions arose. But the moduli space of
$SU(2)$ - instantons, constructed by S. Donladson, has two type of
singularities. The first type comes with reducible solutions
(and we avoid this point taking sufficiently generic metric)
while the second type comes from the conformal invariance
of the instanton theory: the so - called bubbling  instantons
exist due to the conformal invariance. So in the Donaldson
theory one deals with geniunly non compact moduli spaces.
In the Seiberg - Witten theory the situation is much more
convenient for investigations: the theory is not conformally
invariant (so the solutions change when we take another
riemannian metric in the same conformal class as our given
$g$) but the moduli spaces are compact! Due to the fruitful computation
firstly proposed in [2] one gets a bound on norms of solutions
and thus one establishes that the moduli spaces are compact.
The computation is so simple that every person can make it
herself. For this one uses the Weitzenbock formula mentioned above
$$
D^*_a D_a \phi = \na^*_a \na_a \phi + \frac{1}{2} F^+_a \phi +
\frac{s}{4} \phi
$$
where $F_a^+$ is considered as an endomorphism of $W^+$ and $s$ is the
scalar curvature of the fixed riemannian metric $g$. Let
$(a, \phi)$ is a solution of the Seiberg - Witten system, then
we take the pointwise  norm $\vert \phi \vert^2$ and at the point where it
is maximal one gets
$$
\aligned
0 \leq \De \vert \phi \vert^2 =
2<\na_a^* \na_a \phi, \phi> - 2<\na_a \phi, \na_a \phi> \\
\leq 2 <\na^*_a \na_a \phi, \phi> = - \frac{s}{2} \vert \phi \vert^2
\\ - <(\phi \otimes \bar \phi)_0 \phi, \phi> =
- \frac{s}{2} \vert \phi \vert^2 - \half \vert \phi \vert^4 \\
\endaligned
$$
and if $\vert \phi \vert^2$ is non trivial at the maximal point one can
obtain dividing by $\vert \phi \vert^2$  the bound on the norm $\vert \phi
\vert^2$:
$$
max_{x \in X} \vert \phi \vert^2 = max (0, -s)
$$
(see [2]). This means additionally that if $g$ has non negative scalar
curvature then there are no irreducible solutions for the Seiberg - Witten
system of equations. Therefore in this case the moduli space for a generic
riemannian metric is empty.

Now we know that for a generic metric the moduli space is a smooth
compact real manifold. But what about  the  orientability? The operators
which represent the linearization of the Seiberg - Witten equations
are $D_a$ and $\de_a$ and the first one defines the canonical orientation
being originally a complex operator. The second can be endowed with an
orientation  by the choice of an orientation of the trivial determinant
line $det H^0 (X, \R) \otimes  det H^1 (X, \R) \otimes H^+(X, \R)$.
This data is the same as in the Donaldson theory so the moduli spaces
of instantons and monopoles can be oriented simultaneously.

The definition of the Seiberg - Witten invariant requires that
the based manifold $X$ has $b_2^+(X) > 1$. First
of all one considers the case when the $Spin^{\C}$ - structure $c$
induces zero dimensional moduli spaces for generic metrics.
Then for the generic metric $g$ the moduli space $\sM(g,c)$
is a compact oriented real manifold of dimension 0 so it is
a finite set of points. One can count the number of these points with
signs defined by the orientation getting an integer number
$N_{SW}(g, c)$. But since $b_2^+(X) > 1$ this number doesn't
depend on the choice of generic riemannian metric: if $g_0, g_1$ are
two generic metric with the numbers $N_{SW}(g_0, c)$
and $N_{SW}(g_1, c)$ which can be jointed by a smooth path $\ga$ (see above)
and the function $N_{SW}(\ga, c)$ has to be continuous
under the deformation along $\ga$. But any integer continuous function
has to be constant. This implies that two integer numbers $N_{SW}(g_0, c)$
and $N_{SW}(g_1, c)$ are the same. It means that if $b_2^+(X)>1$
the number $N_{SW}$ doesn't depend on the choice of generic riemannian
metric and is an invariant of the smooth structure. Thus in this case
we understand $N_{SW}$ as a function on the set of appropriate
$Spin^{\C}$ - structures with integer values. The simple case
with zero dimensional moduli spaces just corresponds to complex geometry
--- we've reproved Wu theorem and we know that the case when
$$
v.dim \sM = 0 = \frac{1}{4}(c^2 - 2 \chi - 3 \si)
$$
is exactly the complex case. Further let us define the function
$N_{SW}$ on the subset of $Spin^{\C}$ - structures with
negative $c_2(W^+)$ just as zero. If for $Spin^{\C}$ structure
$c$ this number is positive and even then there exists a procedure
which again gives us an integer number which we understand as the value
of $N_{SW}$ on this class. Namely the ambient space $\Cal C^*$
of irreducible pairs modulo gauge equivalence always has a distinguished
2 - cohomological class which can be derived as follows. Let us fix
a point $x \in X$ and consider the subgroup $\sG_0$ of the gauge group
$\sG$ consists of all gauge transformations identical on the fiber
over $x$. Consider the quotient space
$$
\Cal C_0^* = (\sA_h(det W^+) \times \Ga(W^-))^* / \sG_0
$$
which is a principle $U(1)$ - bundle over $\Cal C^*$. The canonical action
is induced by $U(1)$ transformations of the fiber over $x$. This
natural principle $U(1)$ - bundle has a topological characteristic ---
the first Chern class $\eta \in H^2(\Cal C, \Z)$. On the other hand
since $v.dim \sM$ is even it means that one has the fundamental
class $[\sM] \in H_{2d} (\Cal C^*, \Z)$ (perhaps after a small deformation)
where
$$
2 d = \frac{1}{4}(c^2 - 2 \chi - 3 \si)
$$
is even. Then one takes the natural pairing
$$
<\eta^d, [\sM]> \in \Z
$$
and it is an integer number which  again depends (if $b_2^+(X) > 1$)
only on smooth structure. So this number is the value of $N_{SW}$
if $\frac{1}{4}(c^2 - 2 \chi - 3 \si)$ is even. In the odd case
we again say that the number is trivial.

Thus the Seiberg - Witten invariant for a smooth compact orientable
real 4 - dimensional manifold $X$ is an integer function
$$
N_{SW}: Spin \to \Z
$$
where $Spin$ is the set of possible $Spin^{\C}$ - structures over $X$.
This function is invariant under the action of $Diff^+X$ and one knows
how it changes under the action of whole $Diff X$.
At the same time there is a natural symmetry namely if one takes
two $Spin^{\C}$ - structures $\pm c$ where $c$ is a fixed class
then one gets a symmetry between the solutions. But since
the question of compatibility of the orientations depends on
numerical characteristics of $X$ one gets
$$
N_{SW}(-c) = (-1)^k N_{SW}(c)
$$
where
$$
k = \frac{1}{4}(\chi(X) + \si(X)).
$$
It's easy to see that if a $Spin^{\C}$ - structure with even
second Chern class of the spinor bundle exists then $k$ belongs to
integer numbers. If this number is even then the Seiberg - Witten function
$N_{SW}$ is symmetric while in the odd case it is skew symmetric.

If $b_2^+(X)$ is equal to 1 then it remains a possibility to define
the invariant using the chamber structure on the space of all riemannian
metrics. Often this variant is quite useful (see [2]).
The space of riemannian metric is divided  by so called walls into
a set of disjoint chambers. Every wall consists of such riemannian metrics
that harmonic 2 - forms representing the $Spin^{\C}$ - structure class
$c$ are anti self dual. But the number derived from the Seiberg - Witten
moduli space doesn't change inside of any chamber; if two riemannian
metrics $g_0$ and $g_1$ lie in the same chamber we can join them by
a smooth path which doesn't intersect the walls therefore the value
of the integer function has to be constant along the path. Moreover
if $c^2 > 0$ then as we've mentioned above there are no reducible solutions
at all for any riemannian metric. Therefore even if $b_2^+ = 1$ in this case
one can define the invariant integer function without any references to
the chamber structure. But for the main interest of the smooth 4 - geometry
which is focused on complex projective plane considered as a real
4- dimensional manifold this remarks gives almost nothing. The point is
that as we've seen from the Kronheimer and Mrowka calculations
there are no irreducible solutions for the Seiberg - Witten system
if the defining riemannian metric has non negative scalar curvature.
But it's the case of the "canonical" kahler metric on $\C \proj^2$
thus despite of the fact that for every $Spin^{\C}$ - structure $c$
over $\C \proj^2$ there are no metrics admitting anti self dual
realization for the class (the intersection form is positive definite
so every nontrivial class has positive square) and consequently
there is no the chamber structure on the space of all metrics
the invariant has to be trivial. Really since there are no walls
it means that one can join every riemannian metric with the canonical
kahler metric but for the last one the integer number equals to zero
hence it is trivial for all riemannian metrics. This type of arguments
ensures that every 4 - dimensional manifold which admits riemannian metric
with non negative scalar curvature has trivial Seiberg - Witten invariant.

On the other hand we have the original case which was the
starting point of the theory. We mean the spin case when one has over $X$
a distinguished $Spin^{\C}$ - structure $c = 0$. This is possible if
and only if  the second Shtiffel - Whitney class of $X$ vanishes.
In this case it's easy to see that  any riemannian metric admits
anti self dual representation  for this class since the trivial
class is represented by zero harmonic form which of course is anti
self dual (as well as self dual). This feature  leads one to apply
 some additional  gears for the definition and calculations of
the invariant.

 First of all what is the difference between nontriviality of
the invariant and the existence of solutions? The point is that even if
the invariant is trivial then solutions may exist. We will illustrate
this fact using some examples below. But if the invariant is nontrivial
then the existence of solutions is stable; for example for any riemannian
metric a solution exists (may be it is a reducible solution but
nevertheless) if the invariant is non trivial. But there is another kind
of deformation for the system, namely on can perturb the second equation
by a sufficiently small self dual form getting the following system
$$
\aligned
D_a (\phi) = 0 \\
F_a^+ = - (\phi \otimes \bar \phi)_0 + \eta\\
\endaligned
$$
where $\eta$ is a self dual form. The point is that using this deformation
we can define the invariant in the cases mentioned above. For the trivial
$Spin^{\C}$ - structure $c = 0$ let us take an appropriate deformation
of the original system with a 2 - form $\eta$ such that this perturbed system
doesn't admit reducible solutions. If our riemannian metric $g$ has non
negative scalar curvature then the adding of an appropriate $\eta$ kills
the vanishing argument for spinor fields above (see [2]).

Now let us come to the examples which we promised to present as the
illustration of the difference.
To deal with first of all let us recall two basic facts which were
established in the starting work [1]. Using these facts Witten
reproved in [1] the Donaldson result which claims that
some algebraic manifolds aren't  diffeomorphic to their topomodels.
So if $X$ is a kahler manifold with $b_2^+(X) > 1$ endowed as usual
with an integrable complex structure $I$ and kahler metric with
kahler form $\om$ then the Seiberg - Witten invariant of canonical class
$K_I$ as $Spin^{\C}$ - structure equals to $\pm 1$. On the other hand
if one takes the connected sum
$$
Y = X_1 \sharp X_2
$$
where $X_1$ and $X_2$ both have non trivial $b_2^+(X_i)$ then
all the invariants vanish for $Y$. These means that almost all
topomodels have trivial Seiberg - Witten invariants.

We'll prove the first statement in a time in much more
wider context following Taubes results about symplectic
manifolds. The author in [1] argues as follows to prove the
second statement. Studying the gluing procedure one
gets that for a glued metric on $Y$ there is a non trivial
circle action on the moduli space of solutions. For example
if the index formula predicts that the moduli space has dimension
zero then it happens that the moduli space is either 1 - dimensional
for glued metrics or is empty for a generic one. Thus for glued
metrics solution exist while the invariant is trivial and
for generic metric there are no solutions at all.
On the other hand it is quite fruitful technique
in the theory to consider some connected sums instead of
a given manifold, but usually one takes the connected sum
with a number of $\overline{\C \proj}^2$s which is the simplest
manifold with trivial $b_2^+$. We will return to these
connected sums at the end of this lecture.

The results established in the theory can be divided into two sets:
the set of so called vanishing results   and the set of "positive"
results. We've presented two from the first set (non negative
scalar curvature and the connected sum) while the second set has been
represented by the fact that every kahler manifold has non trivial
Seiberg - Witten invariant. This result is quite natural and
expected since the theory looks like a complex theory. Really
in the definition of the invariants the first case is when
one has an almost complex structures over $X$ compatible with
$g$. In this case the moduli spaces are zero dimensional and arithmetically
it corresponds to the condition
$$
c^2 = 2 \chi + 3 \si
$$
so our $Spin^{\C}$ - structure $c$ is a canonical class over $X$.
If the based manifold is simply connected then there is a geometric
procedure which relates the invariant coming from positive dimensional
moduli space and zero dimensional moduli space over a connected sum.
As a pattern we use the following remark from [5]. Let $X$ be
a based manifold as usual in the theory with $b_2^+>1$
and let for a $Spin^{\C}$ - structure $c \in H^2(X, \Z)$
and a riemannian metric $g$ one has the corresponding moduli space
$\sM_X(g, c)$. Glue together our $X$ with $\overline{\C \proj}^2$
getting  another smooth manifold
$$
Y = X \sharp \overline{\C \proj}^2
$$
endowed with $Spin^{\C}$ - structure
$$
c + h  \in H^2(Y, \Z) =  H^2(X, \Z) \oplus \Z h,
$$
where $h$ is the standard generator of $H^2(\overline{\C \proj}^2, \Z)$
with square -1 and a riemannian metric which is the gluing of our given
$g$ and the standard metric on reverse oriented projective plane.  For
this $Spin^{\C}$ - structure $\tilde{c}$ and riemannian metric $\tilde{g}$
one has the corresponding moduli space $\sM_Y(\tilde{g}, \tilde{c})$.
First of all let us compute the virtual dimension
$$
\aligned
v. dim \sM_Y (\tilde{g}, \tilde{c}) =  \frac{1}{4}(\tilde{c}^2 -
2 \chi(Y) - 3 \si(Y)) = \\
\frac{1}{4}(c^2 - 1 - 2 (\chi(X) + 1) - 3 (\si(X) - 1)) =
\frac{1}{4}(c^2 - 2 \chi(X) - 3 \si (X))
\\
\endaligned
$$
thus the moduli spaces have the same virtual dimensions. This gives the hint
that they are isomorphic. The comparing of solutions form these moduli spaces
uses some special technique which is well known from the Donaldson theory
(shrinking the neck etc). Namely changing the canforamlly flat metric on
the neck joining two manifolds one can divide any solution from
$\sM_Y$ into three parts: the part which lives over $X$, the part
living over $\overline{C \proj}^2$ and the remainder over the conformally
flat neck. The first part can be deformed to a solution over $X$ itself
while the part over $\overline{C \proj}^2$ has to be reducible with trivial
spinor field (since the standard metric has positive scalar curvature)
but this reducible solution is unique up to gauge transformations.
It remains to check that one can arrange the smoothing
procedure over the neck to glue the classes from the moduli space
$\sM_X$ and this unique reducible solution. Therefore the moduli
spaces are isomorphic and consequently the invariants for $X$ and
$Y$ are the same.

Now let us  modify this construction.
In [10] one finds the following geometrical remark which is very
closed to this from [5]. Let
$c \in H^2(X, \Z)$ be a $Spin^{\C}$ - structure over a simply connected
based manifold $X$ with an even positive virtual dimension
$$
\frac{1}{4} (c^2 - 2 \chi - 3 \si) = 2 d
$$
where $d \in \Bbb N$.  Then the $Spin^{\C}$ - structure can be extended to
a canonical class over connected sum.
$$
Y = X \sharp d \overline{\C \proj}^2.
$$
Really take
$$
K_c = c + \sum_{i=1}^d 3 h_i
$$
where $h_i$ are the generators over each $\overline{\C \proj}^2$
and compute
$$
\aligned
K_c^2 = c^2 - 9 d = 8 d + 2 \chi(X) + 3 \si(X) - 9 d \\
= 2 \chi(X) + 3 \si(X) - d = 2 \chi(Y) + 3 \si(Y) \\
\endaligned
$$
since
$$
\aligned
\chi(Y) = \chi(X) + d\\
\si(Y) = \si(X) - d.\\
\endaligned
$$
Moreover the extension of $c$ to a canonical class exists if and only
if the number $\frac{1}{4}(c^2 - 2 \chi - 3 \si)$ is even.
Why we search this extension on the connected sums of this type?
The point is that one has to take connected sum with manifolds
whose intersection forms are negative definite. In other cases
one kills all the invariants (see the vanishing result above).
To simplify the construction we would like to take simply connected
manifolds  but due to the Donaldson result every simply connected
4 - manifold with definite intersection form is homeomorphic
to d copies of projective planes. Since we don't know
until now are there some other smooth structures on $\C \proj^2$
unlike to the standard one we could not explain why we take
this concrete smooth $\C \proj^2$ but we believe that the famous
smooth Poincare conjecture is true consequently it would be
only one manifold which can be taken as the summand.

Now the point is that for simply connected based manifold
the invariants computed using $\sM_X(g, c)$ over $X$ and
using $\sM_Y(\tilde{g}, K_c)$ over $Y$ coincide. It means
that one can reduce the computation of the invariants
to an almost complex case. It means that the theory is complex
(or almost complex) itself.

Returning to the definition we recall that for a smooth 4 - manifold
$X$ with $b_2^+(X) > 1$ one constructs an integer function
$$
N_{SW}: Spin \to \Z
$$
(where $Spin$ is the set of all possible $Spin^{\C}$ - structures)
which is invariant under the $Diff^+X$ action. Every class $c \in Spin$
with nontrivial image is called basic class. The set of basic classes
is a finite subset $B_X \subset H^2(X, \Z)$ in general case (for the finiteness
see [1]) and this set is an invariant of smooth structure   ones more.
If all the basic classes are canonical for some almost complex structures
then one calls the smooth manifold $X$ as a manifold of simple type.
In the simply connected case one can always reduce the story to
the simple type: one chooses the basic class with maximal virtual dimension
of the corresponding moduli space and then one takes the connected sum of
$X$ with half of the dimension copies of $\overline{\C \proj}^2$.

Thus we will discuss mostly the cases when one has basic canonical classes.
So we are in the framework of hermitian geometry.

The list of the first results which were established in [1]
can be continued exactly in this way. We left the fact that for
every kahler manifold the canonical class is a basic class without
proof because of the following result which is much more wider. In [3]
Taubes proved that the same is true for any symplectic manifold
(with as usual $b_2^+>1$).   First of all let us recall
how these two cases can be related.

Let $X$ be as usual in the theory real 4 - manifold. And let it admit
almost complex structures. Then one can take any pair $(g, J)$ over
$X$ where $g$ is a riemannian metric and $J$ is a compatible almost
complex structure. Together they define automatically
non degenerated 2 - form $\om$ in the usual way. This 2 - form
is self dual with respect to $g$ and to the orientation given
by $J$ and it has bi - type (1,1) with respect to $J$. So
together we get a hermitian triple $(g, J, \om)$ where
each element can be reconstructed from the others. Then
in this terms one has the following hierarchy
$$
AG \subset KG \subset SG \subset GG \subset DG
$$
where $AG$ is algebraic geometry, $KG$ is kahler geometry, $SG$ is
symplectic geometry, $GG$ is gauge geometry so the geometry
of manifolds with non trivial Seiberg - Witten invariants and
the last one of course is the most general differential geometry of
4 - manifolds. The top subject  $AG$ means that one has a hermitian
triple $(g, J, \om)$ over $X$  which satisfies the following conditions:

i) $J$ is integrable;

ii) $\om$ is closed;

iii) the corresponding 2 - cohomology class $[\om]$ belongs to
integer lattice $H^2(X, \Z)$ in $H^2(X, \R)$.

The next subject $KG$ is defined by a hermitian triple which satisfies
the first two conditions i) and ii). Turning to symplectic geometry
one can see that in this case if $X$ has  symplectic form $\om$
then there exists a cone of compatible almost complex structures
thus taking a structure $J$ from the cone and reconstructing
the corresponding riemannian metric one gets a hermitian triple
which satisfies only the second condition. All these compatible
structures are of the same canonical class which is called associated
canonical class of the symplectic structure. These almost
complex structures can be exploited to define invariants of
smooth structure which are the Gromov invariants
counting pseudoholomorphic curves (we'll discuss
this construction below). Thus
every kahler manifold is a symplectic manifold and it was conjucted
long times ago that the reverse implication takes place as well
before Thurston constructed examples of symplectic manifolds
which do not admit integrable complex structures. Now we get the implication
$SG \subset GG$ due to the Taubes result: for symplectic manifold
(with $b^+_2 > 1$) the canonical class of associated almost complex
structures is a basic class with the invariant equals to $\pm 1$.
The proof is so elegant e and simple that we recall it here
briefly.

First of all if $X$ is a symplectic manifold with symplectic 2 - form
$\om$ then for every compatible riemannian metric $g$ (which induces
the corresponding almost complex structure $J$) there is a preffered
$Spin^{\C}$ - structure which is
$$
\aligned
W^+ = I \oplus \La^{0,2}\\
W^- = \La^{0,1} \\
\endaligned
$$
with $c = - K_J$ as is clear from the formula. Here $I$ is trivial complex line
bundle, $\La^{0,l}$ is the bundle on (0, l) - forms with respect to
our almost complex structure $J$. The properties of $Spin^{\C}$ - structure
can be illustrated by exact formula in this case: if $\rho$ is a real
1 - form and $\al \oplus \be$ belongs to $\Ga (W^+)$ then
$$
\rho (\al \oplus \be) = \al \rho^{0,1} \oplus \La (\rho \wedge \be)
$$
where $\La$ is the adjoint operator to the wedge multiplication by $\om$:
$$
\La: \Om^{p,q} \to \Om^{p-1, q-1}.
$$
The summand $\La^+$ in the Hodge decomposition takes the form
$$
(\La^{2,0} \oplus \La^{0,2})_{\R} \oplus \R <\om>;
$$
it acts on $W^+$ as follows. The symplectic form $\om$ acts as the diagonal
operator with eigenspaces $I$ and $\La^{0,2}$ with eigenvalues
$-i$ and $i$ respectively. Each other 2 - form $\rho \in (\Om^{2,0}
oplus \Om^{0,2})_{\R}$ acts as
$$
\rho(\al \oplus \be) = \al \rho^{0,2} \oplus \La(\rho \wedge \be);
$$
from this it's clear that if $X$ is topological K3 surface then
normalized non vanishing 2 - form $\theta \oplus \bar \theta$
from $(\Om^{2,0} \oplus \Om^{0,2})_{\R}$ defines the second complex
structure on $W^+$.

Moreover, the Dirac operator coupled by $a \in \sA(K^{-1})$
has the form
$$
D_a = \partial_a \oplus \bar \partial_a: I \oplus \La^{0,2}
\to \La^{0,1}
$$
so is a coupling of the standard operator $\de = \bar \partial \oplus
\bar \partial^*$. Thus the system in these "coordinates"
reads as
$$
\aligned
\bar \partial_a \al + \bar \partial^* \be = 0\\
F^{0,2}_a = \bar \al \be\\
i F_a^{1,1} \dot \om = (\vert \be \vert^2 - \vert \al \vert^2) \om
\wedge \om. \\
\endaligned
$$
On the other hand (see [3]) there exists such hermitian connection
$a_0 \in \sA_h(det W^+) = \sA_h(K^{-1})$ that the corresponding
full covariant derivative $\na_{a_0}$ on whole $W^+ = I \oplus K^{-1}$
projects to the first summand as the ordinary $d$ and such $a_0$
is unique up to gauge transformations. This is true for absolutely
general almost complex situation. Now if one takes a constant function
and define the following spinor field
$$
\phi = \al \oplus \be = c \oplus 0
$$
then the action of the covariant derivative gives us
$$
\na_{a_0} (c \oplus 0) = 0 \oplus b
$$
where $b$ is a section of $\La^{0,2} \otimes T^* X$. Indeed, the projection
of the resulting section to the first summand $I \otimes T^*X$ has to be
trivial by the definition of $a_0$ hence we get this $b$. This $b$ is
essentially the torsion of our almost complex structure; one can define
the Nijenhuis tensor $N \in Hom (\La^{0,1}, \La^{2,0})$ using this $b$.
It means that this $b$ is trivial if and only if our complex
structure $J$ is integrable. Formally here is the end point of joint
running for the ways of kahler geometry and symplectic geometry
in the theory. But it is only formally because of a hidden additional
argument which was found by Taubes.   Namely despite of the nontriviality
of $\na_{a_0}(c \oplus 0)$ the corresponding coupled Dirac operator
$D_{a_0}$ vanishes on the spinor field $c \oplus 0$ if and only if
the almost kahler form $\om$ is closed so in the framework of
symplectic geometry. This means that despite of the fact that for
a symplectic manifold the torsion $b$ is nontrivial while for a kahler
one it vanishes the coupled Dirac operator vanishes on this preffered
spinor field in symplectic case as well as in kahler case endowing  us with the  preffered  solution for the first equation of the system.
Let us check this fundamental fact directly. For this
one acts by the almost kahler form $\om$ on the preffered spinor field
getting
$$
\om (c \oplus 0) = - i (c \oplus 0).
$$
Differentiating this equality by the covariant derivative $\na_{a_0}$
and then using Clifford multiplication we establish that
$$
\aligned
\na_{a_0}(\om) * (c \oplus 0) + \om (\na_{a_0}(c+ 0)) =  -i \na_{a_0}(c+0)\\
\implies \\
C d^* \om  (c \oplus 0) + i D_{a_0} (c + 0)  = - i D_{a_0} (c+0) \\
\endaligned
$$
where $C$ is a constant and  $d^* \om$ acts as a real 1 - form on
the spinor field and transform $c \oplus 0$ to a section of $W^- =
\La^{0,1}$. Thus we get the following expression
$$
D_{a_0} (c \oplus 0) = \frac{i C}{2} d^* \om (c \oplus 0).
$$
The map from $W^+$ to $W^-$ defined by 1 - form $d^* \om$
is nontrivial of full rank unless the case when $d^* \om$
vanishes. Since $\om$ is self dual than the condition
$$
d^* \om = 0
$$
coincides with the closeness condition
$$
d \om = 0
$$
therefore
$$
D_{a_0}(c+ 0) = 0 \quad \quad {\text if and only if} \quad \quad d \om = 0
$$
so in the symplectic case (see [3]) one has a distinguished solution
for the first equation as well as in the kahler case. To compute the invariant
Taubes uses an appropriate perturbation of the original system
namely the following one
$$
\aligned
D_a(\al \oplus \be) = 0 \\
F^{0,2}_a = \bar \al \be \\
i F_a^{1,1} \dot \om = (\vert \be \vert^2 - \vert \al \vert^2 +
\rho^2) + C \\
\endaligned
$$
where $\rho$ is the perturbation  parameter   and $C$ is a topological
constant which is the slant product of two 2 - cohomological classes:
$$
C = <2 \pi c_1(K^{-1}) \dot [\om], [X]> \in \R.
$$
Direct computation (see f.e. [22]) shows that when $\rho$ comes to be
sufficiently large then only one solution (up to gauge transformations)
survives and this solution is the obvious one
$$
a = a_0, \al = \rho, \be = 0.
$$
This gives the Taubes result.
On the other hand one can establish coming in this way  another
important Taubes result which says that if $K$  is the canonical class
associated with symplectic form $\om$ over a symplectic manifold
with $b^+_2 > 1$ then $K \cdot [\om] \geq 0$ and for any basic class
$k$  one has inequality
$$
\vert k \cdot [\om] \vert \leq K \cdot [\om].
$$
Note that while the first result is quite known in the kahler case
(if $b^+_2 > 1$ it means that $p_g > 0$ and consequently the canonical
bundle has to have non negative pairing with the given polarization)
the second restriction is new.

After this fundamental result was established one could suppose
that symplectic geometry covers gauge geometry. It was disproved
almost immediately in [5]. The construction is again
quite simple. Let us take a symplectic manifold $X$ with $b_2^+(X) > 1$
as usual in the theory and consider the following connected sum
$$
Y = X \sharp N
$$
where $N$ has $b_1(N) = b_2^+(N) = 0$ but non trivial fundamental
group. An example  used in [5] is the direct product of a homological
3 - sphere with non trivial fundamental group and $S^1$ when one kills
then by the surgery the generator of $\pi_1(S^1)$ in the product. Then
one gets a 4 - dimensional smooth manifold with finite non trivial
fundamental group and with trivial Betti numbers except $b_0$ and $b_4$.
Due to the fact proved in [5] and discussed above the invariants
for $X$ and for $Y$ are the same and it means that $N$ has nontrivial
invariants. On the other hand it couldn't be symplectic. Really if
it is then its universal covering $\tilde{Y} \to Y$ has to be symplectic
as well as $Y$. Thus $\tilde{Y}$ by the Taubes theorem has to have
non trivial invariants. But  this smooth  4 - manifold
is represented by the connected sum
$$
\tilde{Y} = X \sharp ... \sharp X
$$
being the universal covering where the number of summands in
this connected sum equals to the order of $\pi_1(N)$.
Therefore the invariants of $\tilde{Y}$ all are trivial since
it is a connected sum of manifolds with positive $b_2^+$s.
Hence one gets the contradiction with the hypoteze that $Y$ is
symplectic (see [5]). It means that the inclusion
$$
SG \subset GG
$$
is proper at least for simply connected based manifolds.

Now we would like to formulate a natural problem
which arises in hermitian or almost kahler 4 - dimensional
geometry. As we've seen above the Seiberg - Witten theory is
complex. It can be reformulated in terms of almost complex
structures and it has a good completely computed
illustration in the framework of symplectic geometry
(we mean the geometry of pseudo holomorphic curves which we
discuss at the end of this lecture).

So the problem is
\proclaim{Problem} Find a geometrical condition on hermitian triple
$(g, J, \om)$ over a smooth almost complex based manifold which
would be equivalent to the nontriviality condition for
the Seiberg - Witten invariant of the canonical class $K_J$.
\endproclaim

Let us repeat again that may be the problem in the simply connected
case has been solved by the Taubes result.

At the same time there were some attempts to find such kind of
constructive condition  in terms of a special map which exists
over each almost complex 4 - dimensional manifold. Recall from [9]
the definition: let $X$ be an almost complex 4 - manifold (so the
manifold which admits almost complex structures) and $\Cal M_X$
is the space of all possible hermitian triples over $X$. Then
one has a map
$$
\tau: \sM_X \to H^2(X, \R)
$$
which is defined absolutely canonically. Namely if $(g, J, \om)$
is a triple belongs to $\sM_X$ then  the first element $g$ and the second
$J$ defines a Hodge star operator (one takes the conformal class of
$g$ and the orientation coming with $J$) hence one can decompose
the third element with respect to $*_{(g, J)}$ using the Hodge
result
$$
\om = \om_H + d \rho_1 + d^* \rho_2.
$$
The first summand $\om_H$ is a harmonic form (moreover due to
self duality of $\om$ it is self dual harmonic form) hence it
represents a 2 - cohomological class $[\om_H] \in H^2(X, \R)$.
This gives us
$$
\tau(g, J, \om) = [\om] \in H^2(X, \R)
$$
and totally the map
$$
\tau: \sM_X \to H^2(X, \R).
$$
This map is called canonical map because it doesn't require
for the definition any additional choices. This naturality implies
the fact that this map is equivariant with respect to the action
of $Diff X$.

If $X$ is algebraic or kahler or symplectic then it's clear that
there exists a triple over $X$ with nontrivial image in $H^2(X, \R)$
namely if $\om$ is closed then
$$
\om = \om_H
$$
and the image $\tau (g, J, \om)$ is exactly $[\om] \in H^2(X, \R)$.
On the other hand it was shown in [9] that the manifold $Y$
constructed in [5] which doesn't admit symplectic structure
at the same time admits a hermitian triple with non trivial image
in $H^2(X, \R)$ (moreover this class is exactly the image
of $[\om]$ with respect to the isomorphism
$$
H^2(X, \Z) \to H^2(Y, \Z),
$$
where $\om$ is the symplectic form on $X$, see above). The next step
was done in [11] where one decoded the Seiberg - Witten system in
spirit of this canonical construction in almost complex geometry.
Namely one has the fact that if riemannian metric $g$ doesn't admit
reducible solution for the Seiberg - Witten system then every solution
(irreducible of course) gives us a hermitian triple with nontrivial
image with respect to $\tau$. One have some additional results in this
way (see [11]).

On the other hand dealing with this canonical map $\tau$ one can generalize
Kahler - Einstein equation which is equivalent to Einstein equation
for Kahler metrics (see [13]). Namely let us find an appropriate
$Diff^+ X$ - invariant equation on the space of hermitian triples
over an almost complex 4 - dimensional $X$. Such an equation
could be interesting in absolutely general setup of classical gravity
theory. After a momoment we see that one has an equation which has the form
$$
\tau (g, J, \om) = K_J
\tag **
$$
where $K_J$ is the canonical class of $J$. Equivariance of $\tau$
with respect to $Diff^+ X$ gives us that this equation is invariant under
the action. On the other hand if $X$ is a Kahler manifold and
$(g, J, \om)$ is the Kahler triple then this equation coincides
(up to constant)  with the Kahler - Einstein equation.
Let us note that a priori solutions to (**) exist if $K_J^2 \geq 0$
which gives us well known Hitchin - Thorpe inequality
for Einstein manifolds (see [15]). So if we consider
a gauge theory with $Diff^+ X$ as the gauge group using the equation
(**) then we get the following result (see [11]) relating this
gauge theory with a "classical"  gauge theory which is the Seiberg - Witten
theory. Namely solutions to (**) exist if $K_J$ is a basic canonical class
over $X$ and $K_J$ can be realized by a self dual harmonic form with respect
to a riemannian metric $g$.
One expects that pure cohomological equation (**) is related to some
"real" equation on forms, that the last one could be related with 
the classical Einstein equation and together these would imply
some new facts about Einstein metrics in dimension 4.

We shall recall here two other constructions which relate
the theory with geometrical topics.

First, according to results of R. Fintushel and R. Stern
the Seiberg - Witten theory is closely related to the knot
theory. The relationship is presented in [7]; there one defines
a special polynomial expression on basic classes weighted by
the Seiberg - Witten numbers (values of the integer valued
Seiberg - Witten function)  and the relationship is that
if one changes the based manifold applying surgery along
a knot $\ga$ then the Seiberg - Witten polynomial is multiplied by
the Alexander polynomial of the knot.

Second, in the symplectic case there is a beautiful correspondence between
Seiberg - Witten and Gromov invariants. Let us recall
the construction of the last one. For a symplectic manifold $X$ with
symplectic form $\om$ one fixes sufficiently generic almost
complex structure $J$ which is compatible with $\om$. Then for
any 2 - homological class $\eta \in H_2(X, \Z)$  one takes pseudoholomorphic
realizations of this class that is such smooth submanifolds
$\Si \subset X$ representing $\eta$ that at each point $s \in \Si$
the tangent space $T_s\Si$ is a complex subspace with respect to
$J$ in $T_sX$. Such a submanifold was called by M. Gromov "pseudoholomorphic
curve" since  in the case when $X, J$ is a complex surface this definition
gives just holomorphic curves which lie on the surface. Thus for
a fixed $J$ one defines a set
$$
\Cal H_{\eta} = \{ \Si \subset X | [\Si] = \eta \in H_2(X, \Z)
{\text and} \forall s \in \Si J(T_s \Si) = T_s \Si \}
$$
and one can prove (see f.e. [6]) that for sufficiently generic
$J$ this set is a real analityc orientable manifold of real dimension
$$
\mu(\mu - K_J)
$$
where $\mu \in H^2(X, \Z)$ is Poincare dual to $\eta$ cohomological
class. One can compare this formula with the classical Riemann - Roch
formula for divisors on complex surface. Comparing these two
formula one has to remember that the first one is real dimension
while Riemann - Roch gives complex dimension of the corresponding
say complete linear system. The difference however is quite
striking: one can see that canonical divisor has in symplectic setup
unique representation by pseudoholomorphic curves while in algebraic
geometry on has the number $p_g$ which can be different from 1.
To define the Gromov invariant one takes a set of
$\half \mu (\mu - K_J)$ points in general position and
derive from $\Cal H_{\eta}$ subspace consists of such pseudoholomorphic
curves which contain all of these chosen points. This subspace consists
of isolated points which can be counted with signs with respect
to a natural orientation and this counting gives us an integer number
which is the Gromov invariant of the given symplectic manifold.
Again we extend this definition varying homological class
$\eta$ inside of the lattice $H_2(X, \Z)$ which gives us
an integer valued function on the lattice which is invariant
under $Diff X$ - action.

The definition of this invariants is very geometrical and
illustrative. Here we deal with natural geometrical objects
which take place in almost complex geometry as well as in symplectic
geometry. Really for any hermitian triple $(g, J, \om)$ one can take
the set of pseudoholomorphic curves representing a fixed
homological class but in this general case there are no results
which would support the extension of the definition. Returning to
the symplectic case   we recall that Taubes proved that
in this case the Seiberg - Witten invariants are equivalent to
the Gromov invariants. As it was mentioned above for symplectic
manifold $X, \om$ there is a preffered $Spin^{\C}$ - structure
$c = - K_{\om}$ (see above); then each other $Spin^{\C}$ - structures
(we suppose that $H^2(X, \Z)$ has no 2 - torsion) can be labeled
by  2 - homological  classes namely  on the level of
spinor bundles we have relationship
$$
W_{\eta}^{\pm} = L_{\eta} \otimes W_0^{\pm}
$$
where $L_{\eta}$ is the line bundle whose  first Chern class is Poincare
dual to homological class $\eta$. It's clear that homological
class $\eta$ corresponds to $Spin^{\C}$ - structure $2 \eta^* - K_{\om}$.
Then the Taubes result looks like the exact correspondence
$$
N_{SW}(2 \eta^* - K_{\om}) = \pm Gr (\eta)
$$
where $Gr$ is the integer Gromov invariant of $\eta$.

Thus while the result of Fintashel and Stern relates the 4 - dimensional
invariants with the geometry of 1 - dimensional submanifolds the Taubes
correspondence $SW = Gr$ defines a link with the geometry of 2 - dimensional
submanifolds.

The correspondence $SW = Gr$ can be exploited in the studying of the last
problem which we mention in these lectures. This problem was stated
and partially solved by Taubes (see [20]). This is the problem of
non symplectic blow down of symplectic manifolds.

Let $X$ be a symplectic manifold. Then if we take the connected sum
$$
Y = X \sharp \overline{\C \proj}^2
$$
then it is symplectic: the construction was called
"symplectic blow up procedure", it corresponds exactly
to algebro geometric blown up procedure. Moreover Taubes
proved that if $Y$ is a symplectic manifold which contains
a smoothly embedded  2 - sphere with self intersection $-1$
then there exists the following decomposition
$$
Y = X \sharp \overline{C \proj}^2
$$
and the procedure which gives us this $X$ can be called
symplectic blow down procedure (being analogously
to well known process in algebraic geometry). The
conjecture which was proposed by Taubes says that
\proclaim{Problem} For each symplectic 4 - manifold
$Y$ there are no non symplectic blow downs so
$Y$ can be decomposed as $X \sharp N$ where $b_2^+(N) = 0$
if and only if $Y$ admits smoothly embedded 2 - sphere with
self intersection $-1$ and $N$ is  $\overline{\C \proj}^2$.
\endproclaim

This problem has many colors and reflects almost all the problems
listed above. We end this notes with the problem which is
from our view point quite important and intrinsic for
the theory itself. As we mentioned above in the first
lecture one expects that a "twistor space"
for the Seiberg - Witten theory can be constructed
in symplectic category and it implies the interest
to each fact related to symplectic geometry.
On the other hand one can see that all fundamental
questions in 4 - geometry form some background of the problem.
For example why one takes these $\overline{\C \proj}^2$s as summands
glued to a given manifold? The reason is that every simply connected
smooth 4 - manifold with negative definite intersection form
is homeomorphic to $k \overline{\C \proj}^2$ and if one
proves that there exists unique smooth structure on $\C \proj^2$
it would imply that every such manifold is not only homeomorphic
to this connected sum but is diffeomorphic to the last one.
This would imply that there is only one operation in the simply
connected case which doesn't kill the invariants. Thus this
problem is closely related to the uniqueness  problem for
the known 4 - manifold. Other things were studied in [10]:
let $Y$ be a symplectic manifold which can be decomposed in connected
sum
$$
X \overline{\C \proj}^2
$$
such that the associated canonical class $K_Y$ splits with respect
to the usual representation
$$
H^2(Y, \Z) = H^2(X, \Z) \oplus \Z h
$$
where $h$ is the generator of 2 - cohomology group
of $\overline{\C \proj}^2$ as follows
$$
K_Y = c_X + k h
$$
where $c_X$ is a $Spin^{\C}$ - structure over $X$ with positive dimensional
moduli space $\sM_X$. Then it can take place if and only if either $c_X$
is a canonical class (but then $\sM_X$ has zero virtual dimension)
or
$$
v . dim \sM_X = 2
$$
and $k = 3$. So arithmetically this situation (when one has non symplectic
blow down) can happen. But operating with moduli spaces of Seiberg - Witten
solutions and with moduli spaces of Gromov pseudoholomorphic curves
one gets the contradiction with the condition that  $Y$ is a symplectic
manifold (see [10]).  Thus again the problem is projected on the question
is there some other smooth manifold which has the same topological
properties as $\overline{\C \proj}^2$ but admits different smooth structure?
$$
$$
$$
$$

Let us say again that these notes do not pretend to be a complete survey of
the Seiberg - Witten theory at all. We just remark here some
points of the theory which have clear geometrical description or
definition. And we try to illuminate a geometrical point of view
on the theory at all. We skip a lot of technical details and
quite important results since they can be found in the literature.
One can divide the references to the subject into three groups.
The first one   contains original papers
in which all the most important results were established. 
 The second part ([16] - [21]) consists of textbooks
and courses on the Seiberg - Witten theory. The number of such
courses increases all the time so we can list below 
just the top of the aciberg. The third part is exactly the survey
[22] which has to be recommended for every reader. But this paper
itself requires something to be considered so at the same part we add
a textbook on 4 - dimensional geometry contains all what one needs to
study  the subject in the widest context.


\Refs\nofrills{\bf References}
\widestnumber\key{1000}

\ref
\key 1
\by E. Witten
\paper Monopoles and four - manifolds 
\paperinfo Math. Res. Letters, 1 (1994), 
\pages 769 - 789
\endref

\ref
\key 2
\by  P. Kronheimer, T. Mrowka
\paper The genus of embedded surface in the projective plane
\paperinfo Math. Res. Letters, 1 (1994)
\pages  797 - 808
\endref

\ref
\key 3
\by  C. Taubes
\paper  The Seiberg - Witten invariant and symplectic forms 
\paperinfo  Math. Res. Letters, 1 (1994)
\pages 809 - 822
\endref


\ref
\key 4
\by C. Taubes
\paper More constraints on symplectic forms from Seiberg - Witten invariant
\paperinfo  Math. Res. Letters, 2 (1995)
\pages 9 - 13
\endref


\ref
\key 5 
\by  D. Kotschick, J. Morgan, C. Taubes
\paper Four manifolds without symplectic structure but with nontrivial Seiberg - Witten invariant
\paperinfo Math. Res. Letters, 2 (1995)
\pages 119 - 124
\endref


\ref
\key 6
\by C. Taubes
\paper $GR = SW$: counting curves and connections
\paperinfo J. Differ. Geom., 52 (1999), N. 3
\pages 453 - 609
\endref


\ref
\key 7
\by R. Fintushel, R. Stern
\paper Knots, links, and 4 - manifolds
\paperinfo Invent. Math. 134 (1998) N. 2,
\pages 363 - 400
\endref



\ref
\key 8
\by R. Fintushel, R. Stern
\paper Nondiffeomorphic symplectic 4 - manifolds with the same Seiberg - Witten invariants
\paperinfo Geom. Topol. Monogr., 2, Coventry 1999
\pages
\endref



\ref
\key 9
\by N. Tyurin
\paper Abelian monopoles and complex geometry
\paperinfo Math. Notes 65, 3-4 (1999)
\pages 351 - 357
\endref

\ref
\key 10
\by N. Tyurin
\paper Abelian monopoles: the case of positive dimensional
moduli space
\paperinfo Izvestiya: Math., 64: 1 (2000)
\pages 193 - 206
\endref

\ref
\key 11
\by N. Tyurin
\paper Spaces of Hermitian triples and the Seiberg - Witten equations
\paperinfo Izvestiya: Math., 65: 1
\pages 181 - 205
\endref


\ref
\key 12
\by  M. Atiyah, N. Hitchin. I. Singer
\paper Self - duality in four dimensional riemannian geometry
\paperinfo Proc. Roy.Soc., London, Ser. A, 362 (1978), n. 1711
\pages 425 - 461
\endref


\ref
\key 13
\by P. Griffits, J. Harris
\paper Principles of algebraic geometry
\paperinfo Wiley, NY, (1978)
\pages
\endref


\ref
\key 14
\by Yu. Manin
\paper Gauge field theory and complex geometry
\paperinfo Second edition, Springer - Verlag, Berlin 1997
\pages 
\endref

\ref
\key 15
\by N. Hitchin
\paper On compact 4 - dimensional Einstein manifolds
\paperinfo Jour. Diff.Geom., 9 (1974)
\pages 435 - 442
\endref

\ref
\key 16
\by J. Moore
\paper Lectures on Seiberg - Witten invariants
\paperinfo Lecture Notes in Math. 1629, Springer, Berlin (1996)
\pages 
\endref

\ref
\key 17
\by S. Akbulut
\paper Lectures on Seiberg - Witten invariants
\paperinfo Turkish J. Math. 20 (1996) n.1
\pages  95 - 118
\endref

\ref
\key 18
\by J. Morgan
\paper The Seiberg - Witten equations and applications to the topology of smooth four manifolds
\paperinfo Mathematical Notes, 44, University press, Princeton 1996
\pages
\endref

\ref
\key 19
\by M. Marcolli
\paper Seiberg - Witten  gauge theory
\paperinfo Texts and reading in Math. 17, New Delhi (1999)
\pages 
\endref


\ref
\key 20
\by C. Taubes
\paper Seiberg - Witten and Gromov invariants for symplectic
4 - manifolds
\paperinfo First Intern. Press Lecture ser. 2, Sommerville (2000)
\pages
\endref

\ref
\key 21
\by N. Tyurin
\paper Instantons and monopoles (Russian)
\paperinfo Uspekhi Math. Nauk 57 n. 2 (2002)
\pages 85 - 138
\endref


\ref
\key 22
\by S. Donaldson
\paper The Seiberg - Witten equation and 4 - manifold topology
\paperinfo Bull. Amer. Math. Soc. 33:1 (1996)
\pages 45 - 70
\endref



\ref
\key 23
\by S. Donaldson, P. Kronheimer
\paper The geometry of 4 - manifolds
\paperinfo Oxford Univ. Press, Oxford 1990
\pages 
\endref


\endRefs
\enddocument